\documentclass[fleqn]{mat01}
\usepackage{times,mathtimy,amssymb,latexsym}
\begin{document}

\def\N{\ensuremath{\mathbb{N}}}
\def\Z{\ensuremath{\mathbb{Z}}}
\def\D{\ensuremath{\mathbb{D}}}
\def\I{\ensuremath{\mathbb{I}}}
\def\R{\ensuremath{\mathbb{R}}}
\def\C{\ensuremath{\mathbb{C}}}
\def\T{\ensuremath{\mathbb{T}}}

\def\examp{\trivlist \item[\hskip \labelsep{\it Example.}]}

\newtheorem{defi}{\rm DEFINITION}[ssection]
\newtheorem{corol}{\rm COROLLARY}[ssection]

\setcounter{page}{317} \firstpage{317}

\markboth{A Fern\'andez Valles}{Local quasinilpotence and common
invariant subspaces}

\title{Joint local quasinilpotence and common invariant subspaces}

\author{A FERN\'{A}NDEZ VALLES}

\address{Department of Mathematics,
University of C\'adiz, Avda. de la Universidad, s/n~11402-Jerez de
la Frontera, Spain\\
\noindent E-mail: aurora.fernandez@uca.es}

\volume{116}

\mon{August}

\parts{3}

\pubyear{2006}

\Date{MS received 27 November 2005; revised 3 February 2006}

\begin{abstract}
In this article we obtain some positive results about the
existence of a common nontrivial invariant subspace for $N$-tuples
of not necessarily commuting operators on Banach spaces with a
Schauder basis. The concept of joint quasinilpotence plays a basic
role. Our results complement recent work by Kosiek \cite{6} and
Ptak \cite{8}.
\end{abstract}


\maketitle

\section{Introduction}

Let $T$ be a continuous linear operator defined on a separable
Banach space $X$. Let us say that $T$ is cyclic if $x\in X$ such
that
\begin{equation*}
\hbox{Linear Span}\{T^n x\hbox{\rm :}\ x\in X\}
\end{equation*}
is dense in $X$.

On the other hand, we said that $T$ is locally quasinilpotent at
$x\in X$ if
\begin{equation*}
\lim_{n\rightarrow \infty} \|T^n x\|^{{1}/{n}}=0.
\end{equation*}
The notion of local quasinilpotence was introduced in \cite{1} to
obtain nontrivial invariant subspaces for positive operators.

Positive results about $N$-tuples of operators with a nontrivial
common invariant subspace have been obtained in \cite{2,4,7,8}. In
this article, we extend the results of Abramovich {\it et~al}
\cite{1} to the case of not-necessarily commuting $N$-tuples of
continuous linear operators.

To extend these results it will be essential to introduce the
notion of joint local quasinilpotence for $N$-tuples of operators
(not necessarily commuting). The results complement the results
obtained in \cite{1}.

The main result of this paper appears in \S3 and is stated as
follows:

\setcounter{section}{3}
\begin{theorem}[\!]
Let $T=(T_1,\dots,T_N)$ be a $N$-tuple of continuous linear
operators. If $T$ is positive with respect to a cone $C$ and
$x_0\in C$ exists such that $T$ is joint locally quasinilpotent at
$x_0${\rm ,} then there exists a common nontrivial invariant
subspace for $T=(T_1,\dots,T_N)$.
\end{theorem}

Moreover, using this theorem we deduce new results on nontrivial
common invariant subspaces for $N$-tuples of operators (see
Theorem~3.2, Corollary~3.3). We will conclude this article with a
section including open problems and further directions.

\setcounter{section}{1}

\section{Joint local quasinilpotence}

It is easy to see that if $(T_1,\dots,T_N)$ are $N$ commuting
operators and they are locally quasinilpotent at $x_0\in X$, then
the compositions $T_{i_1}\dots T_{i_m}; 1\leq i_j\leq N$ for all
$j\in\{1,\dots,m\}$ and all $m\in \N$, are locally quasinilpotent
at $x_0$. In fact the intersection of the sets
\begin{equation*}
Q_{T_i}=\{x\in X, \ \ \hbox{such that} \ \ T_i \ \ \hbox{is
locally quasinilpotent at}\ \ x\},
\end{equation*}
is a common invariant manifold.

However if $T_1,\dots, T_N$ are not commuting, the problem becomes
more complicated.

\begin{examp}
Let $T_1, T_2$ be two operators on the Hilbert space $l_2$ defined
in the following way:
\begin{equation*}
T_1 e_n = \begin{cases}
\begin{array}{ll}
e_{n-1}, &{\rm if}\ \ n\geq 2\\
0, &{\rm if} \ \ n=1\end{array}; &T_2 e_n=\displaystyle\frac{1}{n}
e_{n+1},
\end{cases}
\end{equation*}
where $(e_n)_{n\in \mathbb{N}}$ is the canonical basis of $l_2$.
\end{examp}

The operators $T_1$ and $T_2$ are locally quasinilpotent at $e_k$
for each $k\geq 2$, since $T_1^ne_k=0$ for each $n\geq k$, and
therefore $\lim_{n\rightarrow \infty}\|T_1^ne_k\|^\frac{1}{n}=0$.
On the other hand,
$T_2^ne_k=\frac{1}{k(k+1)\cdots(k+n-1)}e_{n+k}$, hence
\begin{equation*}
\lim_{n\rightarrow\infty}\|T_2^ne_k\|^{{1}/{n}}
=\lim_{n\rightarrow\infty}\left(\frac{1}{k(k+1)\cdots(k+n-1)}\right)^{{1}/{n}}=0
\end{equation*}
and therefore $T_1$ and $T_2$ are locally quasinilpotent at $e_k$
for each $k\geq 2$.

However, $T_1T_2$ and $T_2T_1$ are not locally quasinilpotent at
$e_k$ for each $k\geq 2$. Indeed, since $T_1T_2e_k=\frac{1}{k}e_k$,
\begin{equation*}
\lim_{n\rightarrow \infty}\|(T_1T_2)^ne_k\|^{{1}/{n}} =
\lim_{n\rightarrow \infty}\left\|\left(\frac{1}{k}\right)^n
e_k\right\|^{{1}/{n}}=\frac{1}{k}\neq 0.
\end{equation*}
On the other hand, we know $T_2T_1e_k=\frac{1}{k-1}e_k$ and hence
$\lim_{n\rightarrow
\infty}\|(T_2T_1)^ne_k\|^{{1}/{n}}=\frac{1}{k-1}\neq 0$.

Let $T=(T_1,\dots, T_N)$ be an $N$-tuple of continuous linear
operators on a Banach space $X$ not necessarily commuting. Let us
denote by $T^n$ the collection of all possible products of $n$
elements in $T$.

\begin{defi}$\left.\right.$\vspace{.5pc}

\noindent {\rm Let $T=(T_1,\dots, T_N)$ be an $N$-tuple of
continuous linear operators on a Banach space $X$ not necessarily
commuting. Then, we will say that $T$ is uniform joint locally
quasinilpotent at $x_0\in X$ if
\begin{equation*}
\lim_{n\rightarrow\infty}\max_{S\in T^n}\|Sx_0\|^{{1}/{n}}=0.
\end{equation*}

The notion of uniform joint local quasinilpotence is closely
related with the joint spectral radius defined by Rota and Strang
\cite{9}. We can get more information about spectral theory in
several variables in \cite{7}.

Although the results of this article are formulated under the
hypothesis of uniform joint local quasinilpotence, nevertheless,
sometimes it will be possible to replace the above-mentioned
hypothesis by the following weaker property.}
\end{defi}

\begin{defi}$\left.\right.$\vspace{.5pc}

\noindent {\rm Let $T=(T_1,\dots, T_N)$ be an $N$-tuple of
continuous linear operators on a Banach space $X$. Then we will
say that $T$ is joint locally quasinilpotent at $x_0\in X$ if
\begin{equation*}
\lim_{n\rightarrow \infty} \|T_{i_{1}}T_{i_{2}}\cdots
T_{i_{n}}x_0\|^{{1}/{n}}=0,
\end{equation*}
for all $i_{1},\dots, i_{n}\in \{1,\dots,N\}$.}
\end{defi}

The difference between the concept of uniform joint local
quasinilpotence and joint local quasinilpotence is the lack of
uniformity. Next we see some properties of both concepts.

\setcounter{defin}{2}
\begin{proposition}$\left.\right.$\vspace{.5pc}

\noindent Let $T=(T_1, \dots,T_N)$ be an $N$-tuple of continuous
linear operators on a Banach space $X$ and let us suppose that $T$
is uniform joint locally quasinilpotent at $x_0\in
X\backslash\{0\}$. Then for all polynomial $p$ of $m$
variables{\rm ,} such that $p(0,\dots,0)=0$ we have that
\begin{equation*}
\lim_{n\rightarrow \infty}
\|p(T_{i_1},\dots,T_{i_m})^nx_0\|^{{1}/{n}}=0,
\end{equation*}
where $i_j\in\{1,\dots,N\}; j\in\{1,\dots,m\}${\rm ,} that is{\rm
,} the operator $p(T_{i_1},\dots,T_{i_m})$ is locally
quasinilpotent at $x_{0}$.
\end{proposition}

\begin{proof}
Fix $\varepsilon>0$. Let us suppose that $k\in \N$ is the number
of summands of the polynomial $p$, and let us denote by
$c\in\R_{+}$ the maximum of the modulus of the coefficients of
$p$. Then, since $T=(T_{1},\dots,T_{N})$ is uniform joint locally
quasinilpotent at $x_{0}$, there exists $n_{0}\in \N$ such that
\begin{equation*}
\max_{S\in T^n}\|Sx_{0}\|^{1/n}<\frac{\varepsilon}{c k}
\end{equation*}
for all $n\geq n_{0}$.

Now, taking into account that the polynomial $p$ has no
independent term, for all $n\geq n_{0}$,
\begin{equation*}
\|p(T_{1},\dots,T_{N})^nx_{0}\|^{1/n}\leq (k^n c^n \max_{S\in
T^n}\|Sx_{0}\|)^{1/n}\leq \varepsilon,
\end{equation*}
which proves the desired result. \hfill $\Box$
\end{proof}

\begin{remark}
{\rm In fact the condition on the polynomial $p(\theta)=0$ is a
necessary condition in the proof of Proposition~2.3. Indeed, let
$F$ be the forward shift defined on the sequences space $\ell_{2}$
by $Fe_{n}=\frac{1}{n!}e_{n+1}$, for all $n\geq 1$. It is easy to
see that the operator $I+F$ is not locally quasinilpotent at
$e_{1}$ (where $I$ denotes the identity operator).}
\end{remark}

Definitions~2.1 and 2.2 are the natural extensions of the notion
of local quasinilpotence for $N$-tuples of continuous linear
operators. In fact, let us denote
\begin{equation*}
Q_{T_1\dots T_N}=\{x_0\in X\hbox{\rm :}\ (T_1\dots T_N) \ \hbox{is
joint locally quasinilpotent at} \ x_0\}
\end{equation*}
and let us denote by $\textit{UQ}_{T_1\dots T_N}$ the set of all
uniform joint locally quasinilpotent vectors for
$(T_{1},\dots,T_{N})$. Then, we have the following result.

\begin{proposition}$\left.\right.$\vspace{.5pc}

\noindent Let $T=(T_1,\dots, T_N)$ be an $N$-tuple of continuous
linear operators on a Banach space $X${\rm ,} then the sets
$\textit{UQ}_{T_1\dots T_N}$ and $Q_{T_1\dots T_N}$ are common
invariant manifolds by $T_1,\dots,T_N$.
\end{proposition}

\begin{proof}
Clearly, $x\in Q_{T_1\dots T_N}$ implies that $\lambda x\in
Q_{T_1\dots T_N}$. Now let $x, y \in Q_{T_1\dots T_N}$, and fix
$\varepsilon>0$. Then, there exists some $n_0$ such that
$\|T_{i_1}\dots T_{i_n}x\|<\varepsilon^n$ and $\|T_{i_1}\dots
T_{i_n}y\|<\varepsilon^n$ for each $i_1,\dots,i_n\in \{1,\dots,
N\}$ and each $n\geq n_0$. Therefore,
\begin{align*}
\|T_{i_1}\dots T_{i_n}(x+y)\|^{{1}/{n}}\leq (\|T_{i_1}T_{i_2}\dots
T_{i_n}x\|+\|T_{i_1}T_{i_2}\dots
T_{i_n}y\|)^{{1}/{n}}<2\varepsilon
\end{align*}
for all $n\geq 0$. Therefore $x+y\in Q_{T_1\dots T_N}$ and so
$Q_{T_1\dots T_N}$ is a vector manifold.

Fix $x_0\in Q_{T_1\dots T_N}$ and let $T_{k}x_{0}\in
Q_{T_{1},\dots,T_{N}}$. Then
\begin{equation*}
\lim_{n\rightarrow \infty}\|T_{i_1}\dots
T_{i_n}(T_kx_0)\|^{{1}/{n}}=\lim_{n\rightarrow
\infty}(\|T_{i_1}\dots T_{i_m}T_k
x_0\|^\frac{1}{n+1})^\frac{n+1}{n}=0
\end{equation*}
for each $i_j\in\{1,\dots,N\}; j\in \mathbb{N}$ and for each
$k\in\{1,\dots,N\}$. Therefore $Q_{T_1\dots T_N}$ is a common
invariant manifold for $T_1, \dots, T_N$. Similar proof follows
for the set $\textit{UQ}_{T_1\dots T_N}$, and hence we omit it.
\hfill $\Box$
\end{proof}

The above propositions show that if $Q_{T_1\dots T_N}\neq \{0\}$
($\textit{UQ}_{T_1,\dots,T_N}\neq \{0\}$ respectively) and
$\overline{Q_{T_1\dots T_N}}\neq X$
($\overline{\textit{UQ}_{T_1\dots T_N}}\neq X$ respectively), then
$\overline {Q_{T_1\dots T_N}}$ ($\overline{\textit{UQ}_{T_1\dots
T_N}}$ respectively) is a common nontrivial closed invariant
subspace for $T_1,\dots, T_N$. As far as the invariant subspace
problem is concerned, we need only consider the two extreme cases
$\overline{Q_{T_1\dots T_N}}=X$ and $Q_{T_1\dots T_N}=\{0\}$.

\section{Main results}

Let $X$ be a Banach space with a Schauder basis $(x_n)$. We say
that $T=(T_1,\dots,T_N)$ is positive with respect to the cone
\begin{equation*}
C=\left\{\sum_{j=1}^\infty t_j x_j\hbox{\rm :} \quad t_j\geq
0\right\}
\end{equation*}
if $T_j(C)\subset C$ for all $j\in\{1,\dots,N\}$.

Let us see the main result of this paper.

\setcounter{defin}{0}
\begin{theorem}[\!]
Let $T=(T_1,\dots,T_N)$ be an $N$-tuple of continuous linear
operators on a Banach space with a Schauder basis such that
$T=(T_1,\dots,T_N)$ is positive with respect to the cone $C${\rm
,} and let us suppose that $y_0\in C$ exists such that
$T=(T_1,\dots, T_N)$ is joint locally quasinilpotent at $y_0$.
Then there exists a common nontrivial invariant subspace for $T=
(T_1,\dots,T_N)$.
\end{theorem}

\begin{proof}
Let $(x_n)$ be a Schauder basis of the Banach space $X$ and let
$(f_n)$ be the sequence of coefficient functionals associated with
the basis $(x_n)$.

Assume that $T=(T_1,\dots,T_N)$ is joint locally quasinilpotent at
some $y_0\in C\setminus\{0\}$, i.e., $\lim_{n\rightarrow
\infty}\|T_{i_1}\dots T_{i_n}y_0\|^{1/n}=0$ with $i_j\in\{1,\dots,
N\}$ for all $j\in \mathbb{N}$. Let us suppose that $T_iy_0=0$ for
all $i\in\{1,\dots,N\}$. Then $\bigcap_{i=1}^N \ker(T_i)$ is a
common nontrivial invariant subspace for each $T_1,\dots,T_N$.
Thus, we can suppose that $T_{i_0} y_0\neq 0$ for some
$i_0\in\{1,\dots,N\}$. By an appropriate scaling of $y_0$, we can
assume that $0< x_k\leq y_0$ and $T_{i_0}x_k\neq 0$ for some $k$
and some $i_0\in\{1,\dots,N\}$.

Now let $P\hbox{\rm :}\ X\rightarrow X$ denote the continuous
projection onto the vector subspace generated by $x_k$ defined by
$P(x)=f_k(x)x_k$. Clearly, $0\leq P(x)\leq x$ holds for each $0<
x\in X$. We claim that
\begin{equation}\label{uno}
PT_{i_1}\dots T_{i_m} x_k=0
\end{equation}
for each $m\geq 0$. To see this, fix $m\geq 0$ and let
$PT_{i_1}\dots T_{i_m} x_k=\alpha x_k$ for some nonnegative scalar
$\alpha> 0$. Since $P$ is a positive operator and the composition
of positive operators is a positive operator, it follows that
\begin{equation*}
0< \alpha^n x_k=(PT_{i_1}\dots T_{i_m})^nx_k\leq (T_{i_1}\dots
T_{i_m})^nx_k\leq (T_{i_1}\dots T_{i_m})^ny_0.
\end{equation*}
Let us observe that the following inequality is not true because
the norm is not monotone $\alpha^n \|x_k\|\leq \|(T_{i_1}\dots
T_{i_m})^ny_0\|$. However, if we use the fact that $f_k$ is a
positive linear functional, we conclude that
\begin{equation*}
0< \alpha^n=f_k(\alpha^n x_k)\leq f_k((T_{i_1}\dots
T_{i_m})^ny_0).
\end{equation*}
Consequently, $0< \alpha^n\leq \|f_k\|\|(T_{i_1}\dots
T_{i_m})^ny_0\|$, and so
\begin{equation*}
0<\alpha\leq\|f_k\|^{{1}/{n}}\|(T_{i_1}\dots
T_{i_m})^ny_0\|^{{1}/{n}}.
\end{equation*}
From Definition~2.2 we know $\lim_{n\rightarrow \infty
}\|(T_{i_1}\dots T_{i_m})^ny_0\|^{{1}/{n}}=0$. Thus we deduce that
$\alpha=0$, and condition (\ref{uno}) must be true.

Now let us consider the linear subspace $Y$ of $X$ generated by
the set
\begin{equation*}
\{T_{i_1}\dots T_{i_m}x_k\hbox{\rm :}\ m\in\N;
i_j\in\{1,\dots,N\} \quad \hbox{for all} \ j\in \mathbb{N}\}.
\end{equation*}
Clearly, $Y$ is invariant for each $T_j;  j\in\{1,\dots,N\}$ and,
since $0\neq T_{i_0} x_k\in Y$ for some $i_0\in \{1,\dots, N\}$,
we see that $Y\neq \{0\}$. From (\ref{uno}) it follows that
$f_k(T_{i_1}\dots T_{i_m}x_k)x_k=P(T_{i_1}\dots T_{i_m}x_k)=0$,
therefore $f_k(T_{i_1}\dots T_{i_m}x_k)$ for each $i_1,\dots,i_m$.
This implies that $f_k(y)=0$ for each $y\in Y$, and consequently
$f_k(y)=0$ for all $y\in \bar{Y}$, that is, $\bar{Y}\neq X$. The
latter shows that $\bar{Y}$ is a common nontrivial closed
invariant subspace for the operators $T_j; j\in\{1,\dots, N\}$,
and the proof is complete.\hfill $\Box$
\end{proof}

Let $T_1,\dots,T_N$ be joint locally quasinilpotent operators at
$x_0\in C$. Then we can add arbitrary weights to each matrix
representing the operators $T_1,\dots,T_N$ and still guarantee the
existence of a common nontrivial closed invariant subspace.

\begin{remark}$\left.\right.$
{\rm \begin{enumerate}
\renewcommand\labelenumi{(\alph{enumi})}
\leftskip .1pc
\item First, let us observe that if $(T_{1},\dots,T_{N})$ is
joint locally quasinilpotent at $x_{0}$ it is possible to obtain a
closed invariant subspace $F$ (nontrivial) invariant also for
every positive operator $A$ such that $AT_i=T_iA$ for all
$i\in\{1,\dots, N\}$. Indeed, the above proof follows considering
the closed subspace generated by
\begin{align*}
\hskip -1.25pc &\{AT_{i_1}\dots T_{i_m}x_k\hbox{\rm :}\ m \in \N;
i_j \in \{1,\dots,N\}\\[.4pc]
\hskip -1.25pc &\quad\, \forall\ j\in \mathbb{N}, A\
\hbox{positive}\ AT_{i}=T_{i}A (\forall i)\}.
\end{align*}

\item On the other hand, let us mention that the subspace guaranteed
in the above theorem is in fact an invariant nontrivial
ideal.\vspace{-.5pc}
\end{enumerate}}
\end{remark}

In the following theorem, positivity is with respect to the cone
generated by the Schauder basis of the Banach space.

\begin{theorem}[\!]
Let $X$ be a Banach space with a Schauder basis. Assume that the
matrix $A_k=(a_{ij}^k)$ defines a continuous operator $T_k$ for
all $k\in \{1,\dots,N\}${\rm ,} such that the $N$-tuple
$T=(T_1,\dots, T_N)$ is joint locally quasinilpotent at a nonzero
positive vector. Let $(w_{ij}^k); k\in\{1,\dots,N\}$ be $N$-double
sequences of complex numbers. If the weighted matrices
$B_k=(w_{ij}^ka_{ij}^k); k\in\{1,\dots,N\}$ define continuous
operators $B_k; k\in\{1,\dots,N\}${\rm ,} then $B_1,\dots, B_N$
have a common nontrivial closed invariant subspace.
\end{theorem}

\begin{proof}
Let $(x_n)$ be a Schauder basis of the Banach space $X$, and let
$(f_n)$ be the sequence of functional coefficients associated with
the basis $(x_n)$. Assume that the operators $T_1,\dots, T_N$
satisfy $\lim_{n\rightarrow \infty}\|T_{i_1}\dots
T_{i_n}y_0\|^{{1}/{n}}=0; i_j\in\{1,\dots,N\}; j\in \mathbb{N}$
for some positive nonzero vector $y_0$. An appropriate scaling of
$y_0$ shows that there exists $l$ satisfying $0< x_l\leq y_0$. Let
us suppose that $T_k x_l=0$ for all $k\in\{1,\dots,N\}$, then an
easy argument shows that $B_k x_l=0$ for all $k\in\{1,\dots,N\}$,
and $\bigcap_{k=1}^N \ker(B_k)$ is a nontrivial closed invariant
subspace (here we assume that $B_k\neq 0$ for all
$k\in\{1,\dots,N\}$). Thus, we can suppose that $T_{i_0}x_l\neq 0$
for some $i_0\in\{1,\dots,N\}$.

Now, let us denote by $P\hbox{\rm :}\ X\rightarrow X$, the
positive projection defined by $P(x)=f_l(x)x_l$. Arguing as in the
proof of Theorem~3.1, we can establish that $PT_{i_1}\dots
T_{i_m}x_l=0$ for each $m\in \N$, where $i_j\in\{1,\dots,N\}$ for
all $j\in \mathbb{N}$. In particular, we have $f_l(T_{i_1}\dots
T_{i_m}x_l)=0$ for each $i_1,\dots,i_m$. Consequently, for each
$m\in\N$ and for each positive operator $S\hbox{\rm :}\
X\rightarrow X$ satisfying $0\leq S\leq T_{i_1}\dots T_{i_m}$, we
have
\begin{equation}
0\leq f_l(Sx_l)\leq f_l(T_{i_1}\dots T_{i_m}x_l)=0. \label{dos}
\end{equation}
Next, let us consider the vector subspace $Y$ generated by the set
\begin{equation*}
\{Sx_l\hbox{\rm :}\ \hbox{such that} \ S \ \hbox{satisfies} \
0\leq S\leq T_{i_1}\dots T_{i_m} \ \hbox{for some} \ i_1,\dots,i_m
\}.
\end{equation*}
Clearly $Y$ is invariant for each operator $R_k;
k\in\{1,\dots,N\}$ satisfying $0\leq R_k\leq T_k$. Also, from
(\ref{dos}) it follows that
\begin{equation*}
f_l(y)=0
\end{equation*}
for all $y\in \bar{Y}$, that is, $\bar{Y}\neq X$. The latter shows
that $\bar{Y}$ is a nontrivial closed vector subspace of $X$. Let
$A_{ij}^k; k\in\{1,\dots,N\}$ now be the operators defined by
\begin{equation*}
A_{ij}^k(x_j)=a_{ij}^k x_j \quad \hbox{and} \quad A_{ij}^k(x_m)=0
\quad \hbox{for} \ m\neq j.
\end{equation*}
Since $A_{ij}^k$ satisfy $0\leq A_{ij}^k\leq A_k$ for all
$k\in\{1,\dots,N\}$, it follows that $\bar{Y}$ is invariant for
each one of the operators $A_{ij}^k$. Therefore, the vector
subspace $\bar{Y}$ is invariant under the operators
\begin{equation*}
B_n^k=\sum_{i=1}^n\sum_{j=1}^n w_{ij}^kA_{ij}^k.
\end{equation*}
However, the sequence of operators $(B_n^k); k\in\{1,\dots,N\}$
converges in the strong operator topology to $B_k$. Therefore,
$B_k(\bar{Y})\subset \bar{Y}$ and, thus, the operators $B_1,\dots,
B_N$, have a common nontrivial closed invariant subspace. \hfill
$\Box$
\end{proof}\pagebreak

\setcounter{corol}{3}
\begin{corol}$\left.\right.$\vspace{.5pc}

\noindent Let $X$ be a Banach space with a Schauder basis. Assume
that the positive matrices $A_k=(a_{ij}^k)$ define continuous
operators on $X${\rm ,} which are joint locally quasinilpotent at
a nonzero positive vector. If the continuous operators
$T_k\hbox{\rm :}\ X\rightarrow X$ are defined by the matrices
$T_k=(t_{ij}^k)$ satisfying $t_{ij}^k=0$ whenever $a_{ij}^k=0${\rm
,} then the operators have a common nontrivial closed invariant
subspace.
\end{corol}

\section{Concluding remarks and open problems}

The notion of uniform joint local quasinilpotence is used
extensively in \cite{5} to obtain common nontrivial invariant
subspaces. Both concepts, joint local quasinilpotence and uniform
joint local quasinilpotence, play an important role in the search
of common nontrivial invariant subspaces.

It would be interesting to know something more on the sets
$\overline{Q_{(T_{1},\dots,T_{N})}}$ and
$\overline{\textit{UQ}_{(T_{1},\dots,T_{N})}}$. Our conjecture is
that both sets are equal in majority of the cases.

On the other hand, it would be interesting to extend the results
of Theorems~3.1 and 3.3 for the case of $N$-tuples of positive
operators defined on a Hausdorff topological vector space, where
the partial order is defined by means of a Markushevish basis.

\section*{Acknowledgement}

The authors would like to thank the referee for his comments and
suggestions that have improved this article.

\end{document}